\documentclass{amsart}
\usepackage{tikz-cd}
\usepackage{fancyvrb}
\usepackage[bookmarks=true, linktocpage=true,
bookmarksnumbered=true, breaklinks=true,
pdfstartview=FitH, hyperfigures=false,
plainpages=false, naturalnames=true,
colorlinks=true, pagebackref=true,
pdfpagelabels]{hyperref}
\hypersetup{
	colorlinks,
	citecolor=blue,
	filecolor=blue,
	linkcolor=blue,
	urlcolor=blue
}

\calclayout

\makeatletter
\@namedef{subjclassname@2020}{%
	\textup{2020} Mathematics Subject Classification}
\makeatother

\renewcommand{\S}{\mathrm S}
\newcommand{\comch}{\texttt{ComCH}}
\newcommand{\id}{\mathrm{id}}
\newcommand{\X}{\mathcal X}

\begin{document}
	
\title[A C.A.S. for the study of commutativity up-to-coherent homotopies]{A computer algebra system for the study of commutativity up-to-coherent homotopies}
\author{Anibal M. Medina-Mardones}
\address{Max Plank Institute for Mathematics, Bonn, Germany}
\email{ammedmar@mpim-bonn.mpg.de}
\address{Department of Mathematics, University of Notre Dame, Notre Dame, IN, USA}
\email{amedinam@nd.edu}
\thanks{The author acknowledges financial support from Innosuisse grant \mbox{32875.1 IP-ICT - 1}.}

\begin{abstract}
	The Python package \comch\, is a lightweight specialized computer algebra system that provides models for well known objects, the surjection and Barratt-Eccles operads, parameterizing the product structure of algebras that are commutative in a derived sense. The primary examples of such algebras treated by \comch\, are the cochain complexes of spaces, for which it provides effective constructions of Steenrod cohomology operations at all prime.
\end{abstract} 

\vspace*{-.8cm}

\maketitle

\section{Introduction}

All the basic notions of number, from the integers to the complex, are equipped with a commutative product, and it was believed until Hamilton's introduction of the quaternions, that the product of any number system must be commutative.
Hamilton's discovery encouraged the consideration of other algebraic structures where commutativity was not demanded, and the effect this had on algebra is only comparable to that of non-euclidean geometries on the study of spaces.\footnote{See for example Chapter V of \cite{Kline1982}.} Around a century later, after the development of topology and homotopy, commutativity was revisited and additional levels enriching the basic dichotomy were identified. These correspond to coherent systems correcting homotopically the lack of strict commutativity, and constitute the focus of much current research on theoretical and applied topology.

After the pioneering work of Steenrod \cite{Steenrod47, Steenrod62}, Adem \cite{Adem52}, Serre \cite{Serre53}, Cartan \cite{Cartan55}, Araki-Kudo \cite{ArakiKudo56}, Dyer-Lashof \cite{DyerLashof62}, Stashef \cite{Stasheff63}, Boardman-Vogt \cite{BoardmanVogt73}, May \cite{May70algebraic, May72geometry}, and many others, today there is a rich theory of commutativity up-to-coherent-homotopies whose modern framework is provided by operads and PROPs, and where $E_n$-operads play a central role parameterizing the different levels of homotopical commutativity. In \comch, we focus on the category of chain complexes, and consider two models of the $E_\infty$-operad equipped with filtrations by $E_n$-operads. These are respectively due to McClure-Smith \cite{McClureSmith03} and Berger-Fresse \cite{BergerFresse04} and are known as the surjection and Barratt-Eccles operads.

The homology of algebras over $E_n$-operads are equipped not only with an induced commutative product but also with homology operations when the coefficient ring is the field $\mathbb F_p = \mathbb Z/ p\mathbb Z$.
The study of these operations at the chain level has become an important issue in topological data analysis \cite{medina2018persistence}, condensed matter physics \cite{Kapustin2017}, category theory \cite{medina2020globular} and others areas. To provide researchers with effective tools for their study, \comch\, implements the constructions of \cite{medina2020chain}, making available for the first time chain level representations of these invariants for spaces presented simplicially or cubically. When the prime is $2$, describing Steenrod operations at the chain level is classical and there are implementations for the simplicial \cite{GonReal1999} and cubical \cite{Pilarczyk2016} cases in \texttt{Sage} \cite{sage} and \texttt{ChainCon} \cite{chaincon} respectively. For odd primes, we do not know of any previous implementation either in the simplicial or cubical contexts. In the former case, a different effective approach was developed by Gonzales-Diaz and Real \cite{GonReal2003, GonReal2005} based on the Eilenberg–Zilber contraction.

\section*{Acknowledgment}
We gratefully acknowledge contributions from Djian Post, Wojciech Reise and Michelle Smith. We thank Dennis Sullivan, Kathryn Hess, John Morgan, Greg Brumfiel, Ralph Kaufmann, Paolo Salvatore, Umberto Lupo, Guillaume Tauzin, and Lewis Tunstall for insightful conversations. We also thank the Laboratory for Topology and Neuroscience at EPFL for its support and hospitality while part of this work was developed.

\tableofcontents

\section{Overview of \comch} \label{s: overview}
In this section we describe the overall structure and main functionalities of \comch, referring to its documentation for a description of all of its classes and their methods.\footnote{Currently hosted at \url{https://comch.readthedocs.io/en/latest/}}

\subsection{Free modules and symmetric groups}

Let $R$ be the ring of integers or one of its quotients. In \comch\, the class \texttt{FreeModuleElement} serves to model elements in free $R$-modules, where $R$ is specified by the attribute \texttt{torsion}. Let $\S_r$ be the set of self-bijection of $\{1, \dots, r\}$ regarded as a group by composition. An element $\sigma \in \S_r$ will be represented by the sequence of its values $(\sigma(1), \dots, \sigma(r))$ and it is modeled in \comch\, using the class \texttt{SymmetricGroupElement}.

\subsection{Operads}

Operads parameterize algebraic structures on chain complexes. The precise although lengthy definition can be found for example in \cite{Markl08}. We will present a key example from which the definition can be abstracted. Let $C$ be a chain complex of $R$-module, and consider the set $End^C(r) = Hom(C, C^{\otimes r})$ of $R$-linear maps as a chain complex in the usual way. The set 
\begin{equation*}
End^C = \left\{End^C(r)\right\}_{r \geq 1}
\end{equation*}
is equipped with the following structure: a left action of $\S_r$ on $End^C(r)$ and composition chain maps
\begin{equation*}
\begin{tikzcd}[column sep=small, row sep=tiny]
\circ_i \colon &[-10pt] End^C(r) \otimes End^C(s) \arrow[r] & End^C(r+s-1) \\
& f \otimes g \arrow[r, |->] & (\id \otimes \cdots \otimes g \otimes \cdots \otimes \id) \circ f 
\end{tikzcd}
\end{equation*}
satisfying forms of equivariance, associativity, and unitality.

An $\mathcal O$-coalgebra structure on $C$ is a structure preserving morphism from $\mathcal O$ to $End^C$. We remark that it is also common to consider the operad $End_A$ obtained from the complexes $End_A(r) = Hom(A^{\otimes r}, A)$, referring to operad morphisms $\mathcal O \to End_A$ as $\mathcal O$-algebra structures on $A$. The linear duality functor induces from an $\mathcal O$-coalgebra in $C$ an $\mathcal O$-algebra structure on $A = Hom(C, R)$.

\subsection{Symmetric ring operad}

Let us consider $R[\S] = \left\{R[\S_r]\right\}_{r \geq 1}$ with $R[\S_r]$ the group ring of $\S_r$ thought of as a dg $R$-module concentrated in degree~$0$. It has the structure of an operad with left action induced from left multiplication, and compositions induced from the maps
\begin{equation} \label{eq: compostion of permutations}
\circ_i \colon \S_r \times \S_s \to \S_{r+s-1}
\end{equation}
sending a pair $(x, y)$ to the bijection $x \circ_i y$ represented diagrammatically by
\begin{equation*}
\underbrace{1 \cdots (\overbrace{i \cdots i+s-1}^y) \cdots r+s-1}_x.
\end{equation*}
More precisely, $x \circ_i y$ is the sequence obtained by applying the following three steps: 1) shift up by $s-1$ the values of $x$ greater than $s$, 2) shift up by $i-1$ the values $y$ and, 3) replace the $i$-th value of $x$ with the shifted $y$. We model elements in $R[\S]$ using the class \texttt{SymmetricRingElement} which combines the classes \texttt{FreeModuleElement} and \texttt{SymmetricGroupElement}. For example, we have
\begin{Verbatim}[frame=single, samepage=true]
>>> x = SymmetricRingElement({(2,3,1): -1, (1,3,2): 1})
>>> y = SymmetricRingElement({(1,3,2): 1, (1,2,3): 2})
>>> print(x * y)
- (2,1,3) - 2(2,3,1) + (1,2,3) + 2(1,3,2)
>>> print(x.compose(y, 2))
- (2,4,3,5,1) - 2(2,3,4,5,1) + (1,5,2,4,3) + 2(1,5,2,3,4)
\end{Verbatim}

\subsection{$E_\infty$-operads}

An important class of operads are those defining resolutions of the ground ring $R$ as an $R[\S_r]$-module. Such operads are called \mbox{$E_\infty$-operads}. They typically come equipped with a filtration by so called $E_n$-operads parameterizing different levels of derived commutativity, with $E_1$ corresponding to the lack of any assumed commutativity, and $E_\infty$ to the largest possible degree of homotopical commutativity. \comch\, implements models of two well known $E_\infty$-operads equipped with filtrations by $E_n$-operads which we now describe.

\subsection{Surjection operad}

For a positive integer $r$ let $\mathcal X(r)_d$ be the free $R$-module generated by all functions $x : \{1, \dots, d+r\} \to \{1, \dots, r\}$ modulo the $R$-submodule generated by degenerate functions, i.e., those which are either non-surjective or have a pair of equal consecutive values. There is a left action of $\mathrm S_r$ on $\mathcal X(r)$ which is up to signs defined on basis elements by $\pi \cdot x = \pi \circ x$.
We represent a surjection $x$ as the sequences of its values $\big( x(1), \dots, x(n+r) \big)$. The boundary map in this complex is defined up to signs by
\begin{equation*}
\partial x = \sum_{i = 1}^{r+d} \pm \big( x(1), \dots, \widehat{x(i)}, \dots, x(n+r) \big),
\end{equation*}
and the $i$-th composition $x \circ_i y$ of $x \in \mathcal X(r)$ and $y \in \mathcal X(s)$ is defined, up to signs, as follows. Let $w$ be the cardinality of $x^{-1}(i)$. For every collection of \textit{ordered indices}
\begin{equation} \label{e: order indices}
1 = j_0 \leq j_1 \leq j_2 \leq \cdots \leq j_{w-1} \leq j_w = s
\end{equation}
we construct an associated splitting of $y$
\begin{equation*}
(y(j_0), \dots, y(j_1));\ (y(j_1), \dots, y(j_2));\ \cdots \ ;\ (y(j_{w-1}), \dots, y(j_w)).
\end{equation*}
The element $x \circ_i y \in \X(r+s-1)$ is represented as the sum over the set of order indices \eqref{e: order indices} of the sequence obtained in the following three steps: 1) shift up by $s-1$ the values of $x$ greater than $i$, then 2) shift up by $i-1$ the values of each sequence in the associated splitting of $y$, and finally 3) replace in order the occurrences of $i$ in $x$ by the corresponding sequence in the splitting.

The elements in this operad are modeled using the class \texttt{SurjectionElement}. For example,
\begin{Verbatim}[frame=single, samepage=true]
>>> x = SurjectionElement({(1,2,1,3): 1})
>>> print(x.boundary())
(2,1,3) - (1,2,3)
>>> y = SurjectionElement({(1,2,1): 1})
>>> print({x.compose(y, 1))
(1,3,1,2,1,4) - (1,2,3,2,1,4) - (1,2,1,3,1,4)
\end{Verbatim}

The signs appearing in these constructions are determined by the attribute \texttt{convention} with possible values the strings \texttt{McClure-Smith} and \texttt{Berger-Fresse}. We refer to \cite{McClureSmith03} and \cite{BergerFresse04} for details on these distinct sign conventions.

We will now review the definition of the complexity of a surjection element. The importance of this concept is that the set of surjection elements with complexity less than $n$ defines an $E_n$-suboperad of $\mathcal X$ \cite{McClureSmith03}.

The complexity of a finite binary sequence (i.e. a sequence of two distinct values) is defined as the number of consecutive distinct elements in it. For example, (1,2,2,1) and (1,1,1,2) have complexities 2 and 1 respectively. The complexity of a basis surjection element is defined as the maximum value of the complexities of its binary subsequences. Notice that for elements in $\mathcal X(2)$, complexity and degree agree. The class \texttt{SurjectionElement} models this concept with the attribute \texttt{complexity}. For example,

\begin{Verbatim}[frame=single, samepage=true]
>>> x = SurjectionElement({(1,2,1,3,1): 1})
>>> print(x.complexity)
1
\end{Verbatim}

\subsection{Barratt-Eccles operad}

For a non-negative integer $r$ define the simplicial set $E(\mathrm S_r)$ by
\begin{align*}
E(\mathrm S_r)_n &= \{ (\sigma_0, \dots, \sigma_n)\ |\ \sigma_i \in \mathrm{S}_r\}, \\
d_i(\sigma_0, \dots, \sigma_n) &= (\sigma_0, \dots, \widehat{\sigma}_i, \dots, \sigma_n), \\
s_i(\sigma_0, \dots, \sigma_n) &= (\sigma_0, \dots, \sigma_i, \sigma_i, \dots, \sigma_n).
\end{align*}
It is equipped with a left $\mathrm S_r$-action defined on basis elements by
\begin{equation*}
\sigma (\sigma_0, \dots, \sigma_n) = (\sigma \sigma_0, \dots, \sigma \sigma_n).
\end{equation*}
The chain complex resulting from applying the functor of normalized $R$-chains to it is denoted $\mathcal E(r)$, and the underlying set of the Barratt-Eccles operad is $\mathcal E = \{\mathcal E(r)\}_{r\geq0}$. To define its composition structure we use the Eilenberg-Zilber map. Let us notice that at the level of the simplicial sets $E(\S_r)$ we have compositions
\begin{equation*}
{\circ}_{i}: E(r) \times E(s) \to E(r + s - 1)
\end{equation*}
induced coordinate-wise from $\eqref{eq: compostion of permutations}$.
We define the composition maps on $\mathcal E$ by precomposing
\begin{equation*}
N_\bullet(\circ_i) \colon N_\bullet(E(r) \times E(s))
\longrightarrow
N_\bullet(E(r + s - 1)) = \mathcal E(r+s-1)
\end{equation*}
with the Eilenberg-Zilber map
\begin{equation*}
\mathcal E(r) \otimes \mathcal E(s) =
N_\bullet(E(r)) \otimes N_\bullet(E(s))
\longrightarrow
N_\bullet(E(r) \times E(s)).
\end{equation*}
For example,
\begin{Verbatim}[frame=single, samepage=true]
>>> x = BarrattEcclesElement({((1,2),(2,1)):1, ((2,1),(1,2)):2})
>>> print(x.boundary())
((1,2),) - ((2,1),)
>>> y = BarrattEcclesElement({((2,1,3),):3})
>>> print(x.compose(y, 2))
3((1,3,2,4),(3,2,4,1)) + 6((3,2,4,1),(1,3,2,4))
\end{Verbatim}

The complexity of a Barratt-Eccles element is define analogously to that of surjection elements. In this case too the subset of elements with complexity less than $n$ defines an $E_n$-suboperad \cite{BergerFresse04}.

An important structure present in the Barratt-Eccles operad missing in the surjection operad is a diagonal chain map compatible with compositions. It is given by the Alexander-Whitney diagonal which on basis Barratt-Eccles element is
\begin{equation*}
\Delta(\sigma_0, \dots, \sigma_n) = \sum_{i=1}^n (\sigma_0, \dots, \sigma_i) \otimes (\sigma_i, \dots, \sigma_n).
\end{equation*}
For example,
\begin{Verbatim}[frame=single, samepage=true]
>>> x = BarrattEcclesElement({((1,2), (2,1)): 1})
>>> print(x.diagonal())
(((1, 2),), ((1, 2), (2, 1))) + (((1, 2), (2, 1)), ((2, 1),))
\end{Verbatim}

\section{Steenrod operations} \label{s: steenrod operations}

In this section we effectively describe how to compute Steenrod operations on the mod-$p$ cohomology of spaces presented as simplicial or cubical sets. The implementation of these algorithms is one of the main novel contributions of \comch\, to the available software used in algebraic topology.

\subsection{Steenrod-Adem structures}

Let $\mathrm{C}_r$ be the cyclic group of order $r$ thought of as the subgroup of $\mathrm{S}_r$ generated by an element $\rho$. The elements
\begin{equation*}
T = \rho-1 \quad \text{ and } \quad N = 1+\rho+\cdots+\rho^{r-1}
\end{equation*}
in $R[C_r]$ define a minimal resolution of $R$ by free $R[C_r]$-modules
\begin{equation*}
\mathcal W(r) = R[C_r] \stackrel{T}{\longleftarrow} R[C_r] \stackrel{N}{\longleftarrow} R[C_r] \stackrel{T}{\longleftarrow} \cdots.
\end{equation*}
We denote a preferred basis element of $\mathcal W(r)_i$ by $e_i$.

A Steenrod structure on an operad $\mathcal O$ is a collection, indexed by $r > 0$, of $\mathrm C_r$-equivariant chain maps $\mathcal W(r) \stackrel{\psi_r}{\longrightarrow} \mathcal O(r)$ for which there exists a factorization through an $E_\infty$-operad $\mathcal W(r) \to \mathcal R(r) \to \mathcal O(r)$ such that the first map is a quasi-isomorphism and the second is an $\S_r$-equivariant chain map. If the maps $\mathcal R(r) \to \mathcal O(r)$ define a morphism of operads $\mathcal R \to \mathcal O$, we say the Steenrod structure is a Steenrod-Adem structure.
For any pair of integers $r$ and $i$, a Steenrod structure produces a preferred element $\psi_r(e_i)$ in $\mathcal O(r)_i$.

Steenrod-Adem structures for the surjection and Barratt-Eccles operads are implemented in \comch\, following their introduction in \cite{medina2020chain}. Some examples of $\psi_r(e_i)$ are

\begin{Verbatim}[frame=single, samepage=true]
>>> r, i = 3, 2
>>> y = Surjection.steenrod_adem_structure(r, i)
>>> print(y)
(1,2,3,1,2) + (1,3,1,2,3) + (1,2,3,2,3)
>>> x = BarrattEccles.steenrod_adem_structure(r, i)
>>> print(x)
((1,2,3),(2,3,1),(3,1,2)) + ((1,2,3),(3,1,2),(1,2,3))
\end{Verbatim}

\subsection{Steenrod operations}

Let $A$ be a chain complex of $\mathbb Z$-modules. Let us assume the operad $End_A$ is equipped with a Steenrod structure $\psi \colon \mathcal W \to End_A$. For any prime $p$, define the linear map $D_d \colon A \otimes \mathbb F_p \to A \otimes \mathbb F_p$ by
\begin{equation*}
D_d(a) = \begin{cases}
\psi(e_d)(a^{\otimes p}) & d \geq 0, \\
0 & d < 0.
\end{cases}
\end{equation*}
As in \cite{May70algebraic}, for any integer $s$ the Steenrod operations
\begin{equation*}
P_s : H_\bullet(A; \mathbb F_2) \to H_{\bullet + s}(A; \mathbb F_2)
\end{equation*}
and, for $p > 2$,
\begin{align*}
P_s & \colon H_\bullet(A; \mathbb F_p) \to H_{\bullet + 2s(p-1)}(A; \mathbb F_p), \\
\beta P_s & \colon H_\bullet(A; \mathbb F_p) \to H_{\bullet + 2s(p-1) - 1}(A; \mathbb F_p),
\end{align*}
are defined for a class $[a]$ of degree $q$ respectively by
\begin{equation*}
P_s\big([a]\big) = \big[D_{s-q}(a)\big] \qquad
\end{equation*}
and
\begin{align*}
P_s\big([a]\big) & = \big[(-1)^s \nu(q) D_{(2s-q)(p-1)}(a)\big], \\
\beta P_s\big([a]\big) & = \big[(-1)^s \nu(q)D_{(2s-q)(p-1)-1}(a)\big],
\end{align*}
where $\nu(q) = (-1)^{q(q-1)m/2}(m!)^q$ and $m = (p-1)/2$.

\subsection{Surjections as linear maps}

In this subsection we describe a Steenrod structure on the cochains of simplicial and cubical sets, which defines, by the previous subsection, Steenrod operations on their mod-$p$ cohomology. Using the linear duality functor, it suffices to define a Steenrod structure on chains. Furthermore, by naturality, it suffices to define it on the universal endomorphism operad $End^{C_\bullet}$ with $End^{C_\bullet}(r) = Hom(C_\bullet, C_\bullet^{\otimes r})$, the chain complex of natural transformations from the functor of normalized simplicial or cubical chains to an $r$-iterated tensor product of itself. An element $f$ in this abstract chain complex can be more concretely described as a set $\{f_n\}_{n \geq 0}$ of elements in $C_\bullet(\Delta^\infty)^{\otimes r}$ and $C_\bullet(\mathbb I^\infty)^{\otimes r}$ respectively, with $f_n$ being the element $f\big([0, \dots, n]\big)$ or $f\big([0 ,1]^{\times n} \big)$ respectively. In \comch, elements in $C_\bullet(\Delta^\infty)^{\otimes r}$ and $C_\bullet(\mathbb I^\infty)^{\otimes r}$ are modeled using the \texttt{SimplicialElement} and \texttt{CubicalElement} classes. For example,

\begin{Verbatim}[frame=single, samepage=true]
>>> x = SimplicialElement({((0,1), (1,2,3), (2,3)): 1})
>>> print(x._latex_())
[0,1] \otimes [1,2,3] \otimes [2,3]
>>> y = CubicalElement({((0,1), (2,1), (2,2)): 1})
>>> print(y._latex_())
[0][1] \otimes [01][1] \otimes [01][01]
\end{Verbatim}

To define the Steenrod structure on the universal operad $End^{C_\bullet}$ it suffices to construct a collection of natural $\S_r$-equivariant chains maps $\X(r) \to Hom(C_\bullet, C_\bullet^{\otimes r})$ since we already have a Steenrod-Adem structure on $\mathcal X$. To describe the $\S_r$-equivariant chain map $\X(r) \to Hom(C_\bullet, C_\bullet^{\otimes r})$ we follow \cite{medina2020prop1, medina2018prop2}. Represent a basis surjection element in $\mathcal X(r)_n$ as the labeled directed graph
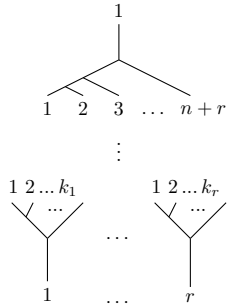
\begin{center}
	\begin{tikzpicture}[scale=.95]
	\draw (0,0)--(0,-.6) node[below, scale=.75]{$1$};
	\draw (0,0)--(.5,.5);
	\draw (-.3, .3)-- (-.2,.5) node[above, scale=.75]{\quad $1\ 2\, ...\, k_1$};
	\draw (-.5,.5)--(0,0);
	\node[scale=.75] at (.11,.4){$...$};
	
	\node[scale=.75] at (1,0){$\cdots$};
	\node[scale=.75] at (1,-.9){$\cdots$};
	
	\draw (2,0)--(2,-.68) node[scale=.75, below]{$r$};
	\draw (2,0)--(2.5,.5);
	\draw (1.7, .3)-- (1.8,.5) node[scale=.75, above]{\quad $1\ 2\, ...\, k_r$};
	\draw (1.5,.5)--(2,0);
	\node[scale=.75] at (2.11,.4){$...$};
	
	\draw (1,2.5)--(1,3) node[scale=.75, above]{$1$};
	\draw (1,2.5)--(0,2) node[scale=.75, below]{$1$};
	\draw (.25,2.125)--(.5,2) node[scale=.75, below]{$2$};
	\draw (.5,2.25)--(1,2) node[scale=.75, below]{$3$};
	\draw (1,2.5)--(2,2) node[scale=.75, below]{\ \quad $n + r$};
	\node[scale=.75] at (1.5,1.75){$\cdots$};
	
	\node[scale=.75] at (1,1.3) {$\vdots$};
	
	\node at (2.85,0){};
	\end{tikzpicture}
\end{center}
where there are no hidden vertices and the strands at the top are joined to the strands at the bottom using the information prescribed by the surjection. Any such graph gives rise to a map in $Hom(C_\bullet, C_\bullet^{\otimes r})$ after associating appropriate maps to the generating pieces
\begin{tikzpicture}[scale=.25]
\draw (0,.5)--(0,1.25);
\draw (0,.5)--(.5,0);
\draw (0,.5)--(-.5,0);
\end{tikzpicture}
and
\begin{tikzpicture}[scale=.25]
\draw (0,.75)--(0,0);
\draw (0,.75)--(.5,1.25);
\draw (0,.75)--(-.5,1.25);
\end{tikzpicture}
in $Hom(C_\bullet, C_\bullet^{\otimes 2})$ and $Hom(C_\bullet^{\otimes 2}, C_\bullet)$ respectively.
Such maps where associated to these generating pieces in \cite{medina2020prop1} for the simplicial case and \cite{medina2020chain} for the cubical one.\footnote{We remark for the interested reader that both of these structures are induced from an $E_\infty$-bialgebra structure on the cellular chains of the interval.} In \comch\, we implement these constructions allowing for surjection elements to act on simplicial and cubical chains. For example, we have
\begin{Verbatim}[frame=single, samepage=true]
>>> x = SurjectionElement({(1,2,1): 1}, convention='McClure-Smith')
>>> a = Simplicial.standard_element(2)
>>> print(x(a))
- ((0,1,2),(0,1)) + ((0,2),(0,1,2)) - ((0,1,2),(1,2))
>>> b = Cubical.standard_element(2)
>>> print(x(b))
- ((2,2),(1,2)) + ((2,1),(2,2)) + ((0,2),(2,2)) - ((2,2),(2,0))
\end{Verbatim}
We remark that in the simplicial context the action of the unique non-degenerate surjection $s_i \colon \{1, \dots, i+2\} \to \{1, 2\}$ with $s_i(1) = 1$ agree up to sign with the cup-$i$ coproducts originally introduced by Steenrod \cite{Steenrod47} and axiomatized in \cite{medina2018axiomatic}.

\subsection{Examples}

We remark that cochains, being defined as $Hom(C_\bullet, R)$, are concentrated in non-positive degrees. We will only give examples in the simplicial context since, in general, the corresponding cubical expressions involve many more terms. We refer to the documentation of \comch\, and its Jupyter notebooks for examples in the cubical context.

1) Let us consider the prime $2$. The value $P_{-1}(x)\big([0,1,2,3,4]\big)$ for $x$ homogeneous of degree $-3$ is the value of $x^{\otimes 2}$ on the following output
\begin{Verbatim}[frame=single, samepage=true]
>>> p, s, q = 2, -1, -3
>>> print(Surjection.steenrod_chain(p, s, q))
((0,1,2,3),(0,1,3,4)) + ((0,2,3,4),(0,1,2,4)) +
((0,1,2,3),(1,2,3,4)) + ((0,1,3,4),(1,2,3,4))
\end{Verbatim}

2) Let us consider the prime $3$. The value $\beta P_{-1}(x)\big([0,1,\dots,8]\big)$ for $x$ homogeneous of degree $-3$ is the value of $x^{\otimes 3}$ on the following output
\begin{Verbatim}[frame=single, samepage=true]
>>> p, s, q = 3, -1, -3
>>> print(Surjection.steenrod_chain(p, s, q, bockstein=True))
2((0,6,7,8),(0,1,2,3),(3,4,5,6)) + ((0,1,7,8),(1,2,3,4),(4,5,6,7))
+ 2((0,1,2,8),(2,3,4,5),(5,6,7,8))
\end{Verbatim}

3) Let us consider the prime $3$ again. The value $P_{-2}(x)\big([0,1,\dots,12]\big)$ for $x$ homogeneous of degree $-4$ is the value of $x^{\otimes 3}$ on the following output
\begin{Verbatim}[frame=single, samepage=true]
>>> p, s, q = 3, -2, -4
>>> print(Surjection.steenrod_chain(p, s, q, bockstein=False))
((0,1,2,3,4),(4,5,6,7,8),(8,9,10,11,12))
\end{Verbatim}

\section{Outlook}

Operations on the homology of algebras that are only $E_n$ for a finite $n$ are well understood homologically \cite{Cohen76} but not at the chain level for $n>2$. The case $n=2$ has been studied by Tourtchine \cite{Tourtchine06} and will be implemented in \comch.

Secondary cohomology operations result from relations among primary. The Cartan and Adem Relations are of particular importance, and constructing cochains enforcing them is open problem for $p > 2$. The case $p = 2$ was treated by Brumfiel, Morgan and the author \cite{medina2020cartan, brumfiel2020cochain} and will soon be implemented in \comch.

Cubical chains appear naturally from simplicial sets through the cobar construction \cite{Baues1980} and in \comch\, we have implemented an $E_\infty$-structure on them. To model the double cobar construction one needs to study permutahedral chains \cite{Kadeishvili2002}. In forthcoming work we describe an $E_\infty$-structure on permutahedral sets suitable for implementation on \comch.

\bibliographystyle{ieeetr} 
\bibliography{bibliography}

\end{document}